\documentstyle{article}
\begin{document}

\input amssym.def
\input amssym.tex

\font\teneusm=eusm10
\newfam\eusmfam
\textfont\eusmfam=\teneusm
\def\scr#1{{\fam\eusmfam\relax#1}}
\def\ad{{\mathrm{ad}}\thinspace}
\newcommand{\adj}{\mbox{\em ad\thinspace}}
\def\l{\left}
\def\r{\right}
\def\la{\langle}
\def\ra{\rangle}

\newtheorem{theorem}{Theorem}[section]
\newtheorem{corollary}[theorem]{Corollary}
\newtheorem{lemma}[theorem]{Lemma}
\newtheorem{proposition}[theorem]{Proposition}
\newtheorem{remark}[theorem]{Remark}

\author{J. Lepowsky and R. L. Wilson}
\title{On Hopf algebras and the elimination theorem for free Lie algebras}
\date{}
\maketitle
\begin{center}
Department of Mathematics, Rutgers University\\
New Brunswick, NJ 08903\\
\end{center}

\begin{center}
e-mail addresses: lepowsky@math.rutgers.edu, rwilson@math.rutgers.edu
\end{center}

\begin{abstract}

The elimination theorem for free Lie algebras, a general principle
which describes the structure of a free Lie algebra in terms of free
Lie subalgebras, has been recently used by E. Jurisich to prove that
R. Borcherds' ``Monster Lie algebra'' has certain large free Lie
subalgebras, illuminating part of Borcherds' proof that the moonshine
module vertex operator algebra obeys the Conway-Norton conjectures.
In the present expository note, we explain how the elimination theorem
has a very simple and natural generalization to, and formulation in
terms of, Hopf algebras.  This fact already follows from general
results contained in unpublished 1972 work, unknown to us when we
wrote this note, of R. Block and P. Leroux.

\end{abstract}

\section{Introduction}

The elimination theorem for free Lie algebras (see \cite{Laz} and
\cite{Bou}) describes how to ``eliminate'' generators; we state and
discuss this theorem in detail below.  It is of interest to us because
it is an important step in E.~Jurisich's proof of her theorem that
R.~Borcherds' ``Monster Lie algebra'' \cite{Bor} is ``almost free'':
This Lie algebra is the direct sum of the 4-dimensional Lie algebra
$\frak g \frak l(2)$ and two $\frak g \frak l(2)$-submodules, each of
which is a free Lie algebra over a specified infinite-dimensional
$\frak g \frak l(2)$-module constructed naturally {}from the moonshine
module (\cite{FLM1}, \cite{FLM2}) for the Fischer-Griess Monster; see
\cite{Jur1}, \cite{Jur2}.  This theorem has been used in \cite{Jur1},
\cite{Jur2} and further, in \cite{JLW}, to illuminate and simplify
Borcherds' proof \cite{Bor} that the moonshine module vertex operator
algebra satisfies the conditions of the Conway-Norton
monstrous-moonshine conjectures \cite{CN}: The McKay-Thompson series
for the the action of the Monster on this structure are modular
functions which agree with the modular functions conjectured in
\cite{CN} to be associated with a graded Monster-module.  In fact,
Jurisich's theorem, including its use of the elimination theorem,
applies to the more general class of generalized Kac-Moody algebras,
also called Borcherds algebras, which have no distinct orthogonal
imaginary simple roots; again see \cite{Jur1}, \cite{Jur2}, and also
\cite{JLW} for further results.

The purpose of this note is to explain how the elimination theorem for
free Lie algebras has a very natural and simple generalization to, and
formulation in terms of, Hopf algebras (and in particular, quantum
groups).

After completing this note we learned from Richard Block that our
results follow from a still more general theorem contained in
unpublished work \cite{BL} of R. Block and P. Leroux, which describes
certain adjoint functors (cf. Remark 3.1 below).  The present note,
then, should be viewed as an exposition showing that the categorical
considerations of Block and Leroux have concrete implications for the
study of the Monster Lie algebra and monstrous moonshine.

Let us write $F(S)$ for the free Lie algebra over a given set $S$ and
$T(S)$ for the free associative algebra (the tensor algebra) over $S$.
The elimination theorem states (see \cite{Bou}) that the free Lie
algebra over the disjoint union of two sets $R$ and $S$ is naturally
isomorphic to the semidirect product of the free Lie algebra $F(S)$
with an ideal consisting of another free Lie algebra.  This ideal is
the ideal of $F(R\bigcup S)$ generated by $R$, and it is the free Lie
algebra $F(T(S)\cdot R)$, where the dot denotes the natural adjoint
action of $T(S)$, which is the universal enveloping algebra of $F(S)$,
in $F(R\bigcup S)$: $$F(R\bigcup S)=F(T(S)\cdot R)
\rtimes F(S);$$ moreover, this $F(S)$-module $T(S)\cdot R$ is in fact
the free $F(S)$-module generated by the set $R$.  (``Module'' means
``left module'' unless ``right module'' is specified.)

The main result of this note is a general structure theorem (which, we
repeat, follows from the results in the unpublished work \cite{BL})
for the algebra freely generated by a given Hopf algebra and a given
vector space, in a sense that we explain below.  (``Algebra'' means
``associative algebra with unit.'')  The theorem expresses this freely
generated algebra as the smash product of the given Hopf algebra with
the tensor algebra over the module for the Hopf algebra freely
generated by the given vector space.  (We also point out that this
freely generated algebra is a Hopf algebra itself in a natural way.)
{}From the special case of this theorem when the given Hopf algebra is
the universal enveloping algebra of some Lie algebra, we obtain as an
immediate corollary a structure theorem expressing the Lie algebra
freely generated by a given Lie algebra and a given vector space as
the semidirect product of the given Lie algebra and the free Lie
algebra over the module for the Lie algebra freely generated by the
given vector space.  Finally, {}from the further special case when
this given Lie algebra is free, we recover the elimination theorem as
an immediate corollary.  (Actually, what we recover is the statement
of the elimination theorem using free Lie algebras over vector spaces
rather than free Lie algebras over sets; the distinction between the
two equivalent formulations of the elimination theorem is simply a
choice of basis for the underlying vector spaces, and such a choice is
not necessary or natural for our purposes.)  It is in this sense that
we are ``explaining'' the elimination theorem as a very special
manifestation of a general Hopf algebra principle.

The first of these two corollaries (the corollary involving a general
Lie algebra) illuminates Jurisich's use of the elimination theorem in
\cite{Jur1} and \cite{Jur2}, and in fact makes her theorem exhibiting
free Lie subalgebras of generalized Kac-Moody algebras, including the
Monster Lie algebra, quite transparent.  For the case of the Monster
Lie algebra, we take the Lie algebra in our setting to be the Lie
subalgebra $\frak g \frak l(2)$ of the Monster Lie algebra, and the
vector space in our setting to be the moonshine module modulo its
(one-dimensional) lowest weight space; then {}from the definition
\cite{Bor}, \cite{Jur1}, \cite{Jur2} of the Monster Lie algebra in
terms of generators and relations, Jurisich's result follows
immediately.  This argument in fact really amounts to her argument in
\cite{Jur2}, with the difference that the Lie subalgebra $\frak g
\frak l(2)$ is viewed as already given, rather than being presented in
terms of generators and relations.  Our main point is that what is
really at work here is a very general, simple and natural Hopf-algebra
principle.  Moreover, our emphasis on tensor algebras (which are the
universal enveloping algebras of free Lie algebras) illuminates the
structure, presented in Section 4 of \cite{JLW}, of certain standard
modules for the Monster Lie algebra as generalized Verma modules---as
tensor algebras over certain modules for the Monster (see especially
Theorem 4.5 of \cite{JLW}).

This work is a continuation of some material in a talk of
J.~Lepowsky's at the June, 1994 Joint Summer Research Conference on
Moonshine, the Monster and Related Topics at Mount Holyoke College, at
which \cite{JLW} was presented.  J.~L. wishes to thank the organizers,
Chongying Dong, Geoff Mason and John McKay, for a stimulating
conference.  The authors are grateful to Elizabeth Jurisich, Wanglai
Li and Siu-Hung Ng, and, as noted above, Richard Block, for helpful
comments related to this work.  Both authors are partially supported
by NSF grant DMS-9401851.

\setcounter{equation}{0}

\section{The results}

In this section, we formulate the main theorem and deduce its
consequences.  We prove the theorem in the next section.

Throughout this note, we work over a field $\Bbb F$.  Let $H$ be a
Hopf algebra, equipped with multiplication $M:H \otimes H
\rightarrow H$, unit $u:\Bbb F \rightarrow H$, comultiplication
(diagonal map) $\Delta:H
\rightarrow H \otimes H$, counit $\epsilon:H \rightarrow \Bbb F$
and antipode $S:H \rightarrow H$, satisfying the usual axioms (cf.
\cite{Swe}).  That is, $H$ is an algebra with multiplication $M$
(typically expressed as usual by juxtaposition of elements) and unit
element $u(1)$ ($1 \in \Bbb F$), which we also write as $1_H$, and is
also a coalgebra such that the coalgebra operations $\Delta$ and
$\epsilon$ are algebra homomorphisms (i.e., $H$ is a bialgebra), and
$H$ is equipped with a linear endomorphism $S$ such that for $h
\in H,$
\begin{equation}
\sum_i h_{1i} S(h_{2i}) = u \epsilon (h) = \epsilon (h) 1_H = \sum_i
S(h_{1i}) h_{2i}.
\label{eq:antipode}
\end{equation}
Here and throughout this note, we use the following conventions for
expressing the comultiplication (in any coalgebra) and its iterates:
For $h \in H$, we write
\begin{equation}
\Delta (h) = \sum_i h_{1i} \otimes h_{2i},
\label{eq:delta}
\end{equation}
where $h_{1i}$ and $h_{2i}$ are suitable elements of $H$; we do not
specify any index set in the summation.  We similarly write
\begin{equation}
(\Delta \otimes 1) \Delta (h) = \sum_i h_{1i} \otimes h_{2i} \otimes
h_{3i},
\label{eq:delta2}
\end{equation}
where the index set and the elements $h_{ji}$ in (\ref{eq:delta}) and
(\ref{eq:delta2}) are unrelated even though the notations are similar.
(Note that by the coassociativity, the expression in (\ref{eq:delta2})
also equals $(1 \otimes \Delta) \Delta (h).$) More generally, for $n >
0$ we write the result of applying the $n$-fold diagonal map to $h$ as
$\sum_i h_{1i} \otimes h_{2i} \otimes \cdots \otimes h_{n+1,i}.$

Note that the defining condition (\ref{eq:antipode}) for the antipode
can also be written:
\begin{equation}
M \circ (1 \otimes S) \circ \Delta = u \epsilon = M \circ (S
\otimes 1) \circ \Delta
\label{eq:antipode2}
\end{equation}
as maps from $H$ to $H$.

The antipode $S$ is an algebra and coalgebra antimorphism (cf.
\cite{Swe}).

Our main object of interest is the algebra freely generated by $H$ and
a given vector space.  First we define the algebra freely generated by
a given algebra and a given vector space.

Given an algebra $A$ and a vector space $V$, let $\scr A(A,V)$
denote the algebra freely generated by the algebra $A$ and the vector
space $V$, in the natural sense specified by the following universal
property:

There is a given algebra homomorphism $A \rightarrow \scr A(A,V)$ and
a given vector space homomorphism $V \rightarrow \scr A(A,V)$ such
that for any algebra $\scr B$ and any algebra homomorphism $A
\rightarrow \scr B$ and vector space homomorphism $V \rightarrow
\scr B$, there is a unique algebra homomorphism $\scr A(A,V)
\rightarrow \scr B$ such that the two obvious diagrams commute.
Such a structure is of course unique up to unique isomorphism if it
exists.

We shall now exhibit such a structure, confirming the existence.  In
the definition, we have not assumed that either $A$ or $V$ is embedded
in $\scr A(A,V)$, but the construction which follows shows that both
of these are in fact embedded in $\scr A(A,V)$.

Set
\begin{equation}
\scr A_n = (A \otimes V)^{\otimes n} \otimes A \mbox{ for } n
\in \Bbb N
\label{eq:An}
\end{equation}
and
\begin{equation}
\scr A(A,V) = \coprod_{n \in \Bbb N} \scr A_n.
\label{eq:A}
\end{equation}
Then $\scr A(A,V)$ is an $\Bbb N$-graded algebra in a natural way; the
product of an element of $\scr A_m$ with an element of $\scr A_n$ is
obtained by juxtaposition followed multiplication of the adjacent
$A$-factors.  This structure clearly satisfies the universal property.
Note that $\scr A(A,V)$ contains $A$ as a subalgebra (in fact, $A =
\scr A_0$) and $V$ as a vector subspace ($V \subset \scr A_1$).  It also
contains the tensor algebras $T(A \otimes V)$ and $T(V \otimes A)$ as
natural subalgebras, but $\scr A(A,V)$ is not the tensor product of
either of these tensor algebras with the algebra $A.$

We shall take $A$ to be the Hopf algebra $H,$ and consider the algebra
$\scr A(H,V),$ which contains $H$ as a subalgebra.

Whenever the Hopf algebra $H$ is a subalgebra of an algebra $C$, we
have the following natural ``adjoint'' action of $H$ on $C$: For $h
\in H$ and $c \in C$,
\begin{equation}
h \cdot c = \sum_i h_{1i}cS(h_{2i}).
\label{eq:adj}
\end{equation}
It is easy to check {}from the definitions and the fact that $S$ is an
algebra antimorphism that
this action makes the algebra $C$ an $H$-module algebra, where $H$ is
now viewed as a bialgebra rather than a Hopf algebra.  This means (cf.
\cite{Swe}) that the linear action $H \otimes C \rightarrow C$
makes $C$ an $H$-module (where $H$ is viewed as an algebra and $C$ as
a vector space) and measures $C$ to $C$ (where $H$ is viewed as a
coalgebra and $C$ as an algebra).  This last condition means that for
$h \in H$ and $c_1, c_2 \in C,$
\begin{equation}
h \cdot c_1 c_2 = \sum_i (h_{1i}
\cdot c_1)(h_{2i} \cdot c_2)
\label{eq:hcc}
\end{equation}
and
\begin{equation}
h \cdot 1_C = u \epsilon (h) =
\epsilon (h) 1_C.
\label{eq:h1}
\end{equation}
Note that while the definition of ``$H$-module algebra'' requires only
that $H$ be a bialgebra and not a Hopf algebra, the definition of the
adjoint action (\ref{eq:adj}) of $H$ on $C$ uses the antipode.  Note
also (the case $C = H$) that under the action (\ref{eq:adj}), any Hopf
algebra is a module algebra for itself.  From the definitions we see
that the product in $C$ and the action (\ref{eq:adj}) are related as
follows: For $h \in H$ and $c \in C$,
\begin{equation}
hc = \sum_i (h_{1i} \cdot c)h_{2i}.
\label{eq:hc}
\end{equation}

We also need the notion of smash product algebra $A \# B$ of an
algebra $A$ which is a $B$-module algebra, $B$ a bialgebra, with the
bialgebra $B$ (cf. \cite{Swe}): As a vector space, $A \# B$ is $A
\otimes B$, and multiplication is defined by:
\begin{equation}
(a_1 \otimes b^1)(a_2 \otimes b^2) = \sum_i a_1(b^1_{1i} \cdot a_2)
\otimes b^1_{2i}b^2
\label{eq:smash}
\end{equation}
for $a_j \in A$, $b^j \in B$.  Suppose that $C$ is an algebra
containing subalgebras $A$ and $H$; that $H$ is a Hopf algebra; that
$A$ is stable under the action (\ref{eq:adj}); and that $C$ is
linearly the tensor product $A \otimes H$ in the sense that the
natural linear map $A \otimes H \rightarrow C$ induced by
multiplication in $C$ is a linear isomorphism.  Then it is clear from
(\ref{eq:hc}) that (\ref{eq:smash}) holds (for $B = H$) and hence that
\begin{equation}
C = A \# H.
\label{eq:CAH}
\end{equation}

For a vector space $V$, we shall write $T(V)$ for the tensor algebra
over $V$.  (Recall that in the Introduction, we also used the notation
$T(S)$ for the free associative algebra over a set $S$; then $T(V)$ is
naturally isomorphic to $T(S)$ for any given basis $S$ of $V$.)

We are now ready to state the main result.  The proof is given in the
next section.  Note that statements (1), (2) and (3) in the Theorem
include three different ``freeness'' assertions.

\begin{theorem}

Let $H$ be a Hopf algebra and $V$ a vector space, and consider the
algebra $\scr A = \scr A(H,V)$ freely generated by $H$ and $V$ and the
canonical adjoint action (\ref{eq:adj}) of $H$ on $\scr A$, making
$\scr A$ an $H$-module algebra.  Then:
\begin{enumerate}

\item  The $H$-submodule $H \cdot V$ of $\scr A$ generated by $V$ is
naturally isomorphic to the free $H$-module over $V$.

\item  The subalgebra of $\scr A$ generated by $H \cdot V$ is naturally
isomorphic to the tensor algebra $T(H \cdot V)$.

\item  As a vector space, the algebra $\scr A$ is naturally isomorphic
to the tensor product
$$\scr A = T(H \cdot V) \otimes H$$
under multiplication in $\scr A$; $T(H \cdot V)$ is an $H$-module
subalgebra of $\scr A$; and with this structure, $\scr A$ is the smash
product of $T(H \cdot V)$ and $H$:
$$\scr A = T(H \cdot V) \# H.$$

\end{enumerate}

\end{theorem}

Here we observe several consequences.  We describe first the special
case in which $H$ is the universal enveloping algebra $U(\frak g)$ of
a Lie algebra $\frak g$ and then the further special case in which
$\frak g$ is the free Lie algebra over a vector space.  For a vector
space $W$, we shall write $F(W)$ for the free Lie algebra over $W$, so
that $F(W)$ is naturally isomorphic to the free Lie algebra $F(S)$ for
any given basis $S$ of $W$.  (As with the notation $T(\cdot)$, we are
using the notation $F(\cdot)$ for both sets and vector spaces.)

\begin{corollary}

Let $\frak g$ be a Lie algebra and $V$ a vector space.  Let $\scr A$
be the algebra freely generated by $\frak g$ and $V$ (defined by the
obvious universal property), or equivalently, the algebra freely
generated by $U(\frak g)$ and $V$, so that $\frak g \subset U(\frak g)
\subset \scr A$.  Let $\frak a$ be the Lie subalgebra of $\scr A$
generated by $\frak g$ and $V$.  Then $\frak a$ is naturally
isomorphic to the Lie algebra freely generated by $\frak g$ and $V$
(in the obvious sense).  Moreover, $\scr A$ is naturally isomorphic to
the universal enveloping algebra of $\frak a$.  Consider the adjoint
action of $\frak g$ on $\scr A$; then the associated canonical action
of $U(\frak g)$ on $\scr A$ is the adjoint action (\ref{eq:adj}).  We
have:
\begin{enumerate}

\item The $\frak g$-module $U(\frak g) \cdot V \subset \frak a \subset
\scr A$ generated by $V$ is naturally isomorphic to the
free $\frak g$-module (i.e., free $U(\frak g)$-module) over $V$.

\item The (associative) subalgebra of $\scr A$ generated by $U(\frak g)
\cdot V$ is naturally isomorphic to the tensor algebra $T(U(\frak g)
\cdot V)$ and in particular, the Lie subalgebra of $\scr A$, and of
$\frak a$, generated by $U(\frak g) \cdot V$ is naturally isomorphic
to the free Lie algebra $F(U(\frak g) \cdot V)$.

\item The algebra $\scr A$ is the smash product
$$\scr A = T(U(\frak g) \cdot V) \# U(\frak g)$$
and in particular,
$$\frak a = F(U(\frak g) \cdot V) \rtimes \frak g;$$
$F(U(\frak g) \cdot V)$ is the ideal of $\frak a$ generated by $V$.

\end{enumerate}

\end{corollary}

\noindent{Proof}:
All of this is clear; to verify the second part of (3), we note that
$F(U(\frak g) \cdot V) \oplus \frak g \subset \frak a$, and the
adjoint action of $\frak g$ on $\scr A$ preserves $F(U(\frak g) \cdot
V)$.  Thus $F(U(\frak g) \cdot V) \oplus \frak g$ is a Lie subalgebra
of $\frak a$ containing $\frak g$ and $V$, giving us the semidirect
product.  $\square$ \vspace{1ex}

\begin{remark}

The Poincar\'e-Birkhoff-Witt theorem has been used here in two places:
It is used to show that $\frak g \subset U(\frak g)$.  It is also used
to show that $\frak a$ is naturally isomorphic to the Lie algebra
freely generated by $\frak g$ and $V$; the proof of the universal
property for $\frak a$---that for any Lie algebra $\frak b$ and any
Lie algebra map $\frak g \rightarrow \frak b$ and any vector space map
$V \rightarrow \frak b$, there exists exactly one Lie algebra map
$\frak a \rightarrow \frak b$ making the two obvious diagrams
commute---uses the fact that $\frak b \subset U(\frak b)$.  Note that
this latter use of the Poincar\'e-Birkhoff-Witt theorem is a
generalization of the similar use of this theorem (in what amounts to
the case $\frak g = 0$) in proving that the Lie subalgebra of $T(V)$
generated by $V$ is the free Lie algebra over $V$ and hence that $T(V)
= U(F(V))$.

\end{remark}

We restate the parts of this corollary involving only the Lie algebra
structure:

\begin{corollary}

Let $\frak g$ be a Lie algebra and $V$ a vector space, and let $\frak
a$ be the Lie algebra freely generated by $\frak g$ and $V$, so that
$\frak g \subset \frak a$ and $V \subset \frak a$.  Then:
\begin{enumerate}

\item The $\frak g$-module $U(\frak g) \cdot V \subset \frak a$
generated by $V$ is naturally isomorphic to the free $\frak g$-module
over $V$.

\item The Lie subalgebra of $\frak a$ generated by $U(\frak g) \cdot V$
is naturally isomorphic to the free Lie algebra $F(U(\frak g) \cdot
V)$.

\item We have the semidirect product decomposition $$\frak a = F(U(\frak
g) \cdot V) \rtimes \frak g,$$ and $F(U(\frak g) \cdot V)$ is the
ideal of $\frak a$ generated by $V$.
$\square$ \vspace{1ex}

\end{enumerate}

\end{corollary}

Now we consider the further special case in which $\frak g = F(W)$,
the free Lie algebra over a given vector space $W$.  We have
immediately:

\begin{corollary}

Let $V$ and $W$ be vector spaces, and consider $F(W), F(V \oplus W)
\subset T(V \oplus W)$.  Then:
\begin{enumerate}

\item The $F(W)$-module
$$U(F(W)) \cdot V = T(W) \cdot V \subset F(V \oplus W)$$
generated by $V$ is naturally isomorphic to the free $F(W)$-module
(i.e., free $T(W)$-module) over $V$.

\item The Lie subalgebra of $F(V \oplus W)$ generated by $T(W) \cdot V$
is naturally isomorphic to the free Lie algebra $F(T(W) \cdot V)$ and
the associative subalgebra of $T(V \oplus W)$ generated by $T(W) \cdot
V$ is naturally isomorphic to the tensor algebra $T(T(W) \cdot V)$.

\item We have the semidirect product decomposition
$$F(V \oplus W) = F(T(W) \cdot V) \rtimes F(W),$$
$F(T(W) \cdot V)$ is the ideal of $F(V \oplus W)$ generated by $V$,
and
$$T(V \oplus W) = T(T(W) \cdot V) \# T(W). \hspace{1em} \square$$

\end{enumerate}

\end{corollary}

\vspace{1ex}

\begin{remark}

The elimination theorem for free Lie algebras over sets rather than
over vector spaces, as recalled in the Introduction, is simply the
obvious restatement of the last corollary for the free Lie algebra
$F(R \bigcup S)$ over the disjoint union of sets $R$ and $S$ in place
of the free Lie algebra $F(V \oplus W)$; we take $R$ and $S$ to be
bases of $V$ and $W$, respectively.  Also, $T(V \oplus W)$ is replaced
by $T(R \bigcup S)$.

\end{remark}

\begin{remark}

Given $H$ and $V$ as in Theorem 2.1, we observe that the algebra $\scr
A = \scr A(H,V)$ freely generated by $H$ and $V$ is naturally a Hopf
algebra.  In fact, define $\Delta:\scr A \rightarrow \scr A \otimes
\scr A$ to be the unique algebra map which agrees with $\Delta:H
\rightarrow H \otimes H$ on $H$ and with the linear map
$V \rightarrow V \otimes V$ taking $v \in V$ to $v \otimes 1 + 1
\otimes v$ on $V$.  Define $\epsilon:\scr A \rightarrow \Bbb F$
to be the unique algebra map which agrees with $\epsilon:H
\rightarrow \Bbb F$ and with the zero map on $V$.  Then $\scr A$
is clearly a bialgebra since $H$ and $V$ generate $\scr A$.  Define
$S:\scr A \rightarrow \scr A$ to be the unique algebra antimorphism
which agrees with the given antipode on $H$ and with the map $-1$ on
$V$.  Then $S$ is an antipode, since (\ref{eq:antipode2}) holds on $H$
and on $V$, and since the set of elements of a given bialgebra on
which (\ref{eq:antipode2}) holds, where $S$ is a given algebra
antimorphism, is closed under multiplication (cf. \cite{Swe}, p. 73).

\end{remark}

\section{Proof of Theorem 2.1}

Now we prove Theorem 2.1.

Consider the free $H$-module $H \otimes V$ generated by $V$ (viewing
$H$ as an algebra) and define the linear map
\begin{eqnarray}
i : H \otimes V & \rightarrow & H \otimes V \otimes H \nonumber \\
    h \otimes v & \mapsto     & \sum_i h_{1i} \otimes v \otimes S(h_{2i})
\label{eq:i}
\end{eqnarray}
($h \in H$, $v \in V$); that is,
\begin{equation}
i = (1 \otimes \scr T) \circ (1 \otimes S \otimes 1) \circ (\Delta
\otimes 1),
\label{eq:i2}
\end{equation}
where $\scr T$ is the twist map:
\begin{eqnarray}
\scr T : K \otimes L & \rightarrow & L \otimes K \nonumber \\
         k \otimes l & \mapsto & l \otimes k
\label{eq:T}
\end{eqnarray}
($k \in K$, $l \in L$) for arbitrary vector spaces $K$ and $L$.
(There should be no confusion between the two uses of the notation
``$i$.'')
Then $i$ is clearly an $H$-module map, where $H \otimes V \otimes H$
is understood as the tensor product of the free $H$-module $H \otimes
V$ and the $H$-module $H$ equipped with the (left) action given by: $h
\cdot k = kS(h)$ for $h, k \in H$.  The image $i(H \otimes V)$
of $i$ is the $H$-submodule of $H \otimes V \otimes H$ generated by
$V$ with respect to this action: $i(H \otimes V) = H \cdot V$.
Furthermore, the $H$-module map
\begin{equation}
i:H \otimes V \rightarrow H \cdot V
\label{eq:i3}
\end{equation}
is an isomorphism (i.e., is an injection), since the linear map
$$1 \otimes 1 \otimes \epsilon : H \otimes V \otimes H \rightarrow H
\otimes V$$
is a left inverse of $i$.  In particular, $H \cdot V$ is a copy of the
free $H$-module generated by $V$.  Recalling (\ref{eq:An}),
(\ref{eq:A}) and writing $\scr A = \scr A(H, V)$ as in the statement
of Theorem 2.1, we observe that $H \cdot V \subset \scr A_1 \subset
\scr A$, and that the first part of Theorem 2.1 is verified.

Now consider the linear map
\begin{eqnarray}
i_1 : H \otimes V \otimes H & \rightarrow & H \otimes V \otimes H
\nonumber \\
  h^1 \otimes v \otimes h^2 & \mapsto & \sum_i h^1_{1i} \otimes v
\otimes S(h^1_{2i})h^2
\label{eq:i1}
\end{eqnarray}
($h^1, h^2 \in H$, $v \in V$), which we may express as the map
\begin{equation}
i_1 : \scr A_1 \rightarrow \scr A_1
\label{eq:i1A}
\end{equation}
given by the composition of $i \otimes 1 : \scr A_1
\rightarrow \scr A_1 \otimes H$ with the multiplication map in
$\scr A$.  We may also write:
\begin{equation}
i_1 = (1 \otimes \scr T) \circ (1 \otimes M \otimes 1) \circ (1
\otimes S \otimes 1 \otimes 1) \circ (\Delta \otimes 1 \otimes 1)
\circ (1 \otimes \scr T).
\label{eq:i1T}
\end{equation}

We now show that $i_1$ is a linear automorphism, and in fact that the
map
\begin{eqnarray}
j_1 : H \otimes V \otimes H & \rightarrow & H \otimes V \otimes H
\nonumber \\
  h^1 \otimes v \otimes h^2 & \mapsto & \sum_i h^1_{1i} \otimes v
\otimes h^1_{2i}h^2
\label{eq:j1}
\end{eqnarray}
is a left and right inverse of $i_1$: For variety, in this argument we
use the symbolism (\ref{eq:antipode2}) and (\ref{eq:i1T}).  Since
\begin{equation}
j_1 = (1 \otimes \scr T) \circ (1 \otimes M \otimes 1) \circ (\Delta
\otimes 1 \otimes 1) \circ (1 \otimes \scr T),
\label{eq:j1T}
\end{equation}
we have
\begin{eqnarray*}
j_1 \circ i_1 & = & (1 \otimes \scr T)(1 \otimes M \otimes 1)
(\Delta \otimes 1 \otimes 1)(1 \otimes M \otimes 1) \cdot\\
& & \hspace{1em} \mbox{} \cdot (1 \otimes S \otimes 1 \otimes 1)(\Delta
\otimes 1 \otimes 1)(1 \otimes \scr T) \\
& = & (1 \otimes \scr T)(1 \otimes M \otimes 1)
(1 \otimes 1 \otimes M \otimes 1) \cdot \\
& & \hspace{1em} \mbox{}\cdot (\Delta \otimes 1 \otimes 1 \otimes 1)(1
\otimes S \otimes 1 \otimes 1)(\Delta \otimes 1 \otimes 1)
(1 \otimes \scr T) \\
& = & (1 \otimes \scr T)(1 \otimes M \otimes 1)
(1 \otimes 1 \otimes M \otimes 1)(1 \otimes 1
\otimes S \otimes 1 \otimes 1) \cdot \\
& & \hspace{1em} \mbox{} \cdot (\Delta \otimes 1 \otimes 1 \otimes 1)
(\Delta \otimes 1 \otimes 1)(1 \otimes \scr T) \\
& = & (1 \otimes \scr T)(1 \otimes M \otimes 1)
(1 \otimes M \otimes 1 \otimes 1)(1 \otimes 1
\otimes S \otimes 1 \otimes 1) \cdot \\
& & \hspace{1em} \mbox{} \cdot (1 \otimes \Delta \otimes 1 \otimes 1)
(\Delta \otimes 1 \otimes 1)(1 \otimes \scr T) \\
& = & (1 \otimes \scr T) (1 \otimes M \otimes 1)
(1 \otimes u \epsilon \otimes 1)
(\Delta \otimes 1 \otimes 1)(1 \otimes \scr T) \\
& = & (1 \otimes \scr T) (1 \otimes \scr T) = 1,
\end{eqnarray*}
where we have used here essentially all the defining
properties of a Hopf algebra: associativity, coassociativity, the
definition (\ref{eq:antipode2}) of the antipode, and the unit and
counit properties.  The fact the $i_1 \circ j_1 = 1$ is verified
similarly.  (The two parts of (\ref{eq:antipode2}) are used in the two
different arguments.)

For $n \geq 1$, define
\begin{equation}
i_n, j_n : \scr A_n \rightarrow \scr A_n
\label{eq:injn}
\end{equation}
by:
\begin{equation}
i_n = (i_1 \otimes 1 \otimes \cdots \otimes 1) \circ \cdots \circ (1
\otimes \cdots \otimes 1 \otimes i_1 \otimes 1 \otimes 1) \circ (1
\otimes \cdots \otimes 1 \otimes i_1),
\label{eq:indef}
\end{equation}
\begin{equation}
j_n = (1 \otimes \cdots \otimes 1 \otimes j_1) \circ (1 \otimes \cdots
\otimes 1 \otimes j_1 \otimes 1 \otimes 1) \circ \cdots \circ
(j_1 \otimes 1 \otimes \cdots \otimes 1).
\label{eq:jndef}
\end{equation}
Using the facts that $j_1 \circ i_1 = 1$, $i_1 \circ j_1 = 1$,
we obtain:
\begin{eqnarray}
j_n \circ i_n & = & 1 \\
i_n \circ j_n & = & 1.
\label{eq:injn=1}
\end{eqnarray}
Note that we may write $i_2:\scr A_2 \rightarrow \scr A_2$,
for example, explicitly as the map
\begin{equation}
h^1 \otimes v_1 \otimes h^2 \otimes v_2 \otimes h^3 \mapsto
\sum_{i,j}\left(h^1_{1i} \otimes v_1 \otimes S(h^1_{2i})\right)
\left(h^2_{1j} \otimes v_2 \otimes S(h^2_{2j})\right)h^3
\label{eq:i2A}
\end{equation}
($h^1, h^2, h^3 \in H, v_1, v_2 \in V$).  Just as the map $i_1$ can be
expressed using the map $i$ and multiplication in $\scr A$ (recall
(\ref{eq:i1A})), the linear automorphism $i_n$ is the composition of
$$i^{\otimes n} \otimes 1:\scr A_n = (H \otimes V)^n \otimes H
\rightarrow (H \otimes V \otimes H)^{\otimes n} \otimes H$$
with the multiplication map in $\scr A$, as in (\ref{eq:i2A}).

We combine the linear automorphisms $i_n$ for $n \geq 0$, where we
take $i_0$ to be the identity map on $H$, to form the graded linear
automorphism
\begin{equation}
\scr I = \coprod_{n \geq 0} i_n : \scr A \rightarrow \scr A,
\label{eq:I}
\end{equation}
so that
\begin{equation}
\left. \scr I \right|_H = i_0 = 1 : H \rightarrow H,
\label{eq:IH}
\end{equation}
the identity map, and
\begin{equation}
\left. \scr I \right|_{H \otimes V} = i : H \otimes V \stackrel{\sim}
{\rightarrow} H \cdot V,
\label{eq:IHV}
\end{equation}
a linear isomorphism and in fact an $H$-module isomorphism (recall
(\ref{eq:i3})).  Moreover, $$\left. \scr I \right|_{T(H \otimes V)} :
T(H \otimes V) \rightarrow
\scr A$$
is a linear isomorphism from the subalgebra $T(H \otimes V)$ of $\scr
A$ to the subalgebra of $\scr A$ generated by $H
\cdot V$, and this subalgebra is naturally isomorphic to the tensor
algebra $T(H \cdot V)$, so that the second statement in the Theorem is
proved:
\begin{equation}
\left. \scr I \right|_{T(H \otimes V)} : T(H \otimes V) \stackrel{\sim}
{\rightarrow} T(H \cdot V);
\label{eq:ITHV}
\end{equation}
in fact, this map is the isomorphism of tensor algebras induced by the
linear isomorphism (\ref{eq:IHV}) of generating vector spaces of the
two tensor algebras.  Thus the linear automorphism $\scr I$ restricts
to algebra isomorphisms (\ref{eq:IH}) and (\ref{eq:ITHV}), which
canonically extends (\ref{eq:IHV}), and the linear decomposition
$$\scr A = T(H \otimes V) \otimes H$$
of the domain of $\scr I$ transports to a linear decomposition
\begin{equation}
\scr A = T(H \cdot V) \otimes H
\label{eq:ATHVH}
\end{equation}
of the codomain; the associated canonical linear map
\begin{equation}
T(H \cdot V) \otimes H \rightarrow \scr A
\label{eq:THA}
\end{equation}
is the map induced by multiplication in $\scr A$.  In particular, the
natural linear map (\ref{eq:THA}) induced by multiplication is a
linear isomorphism.  This proves the first part of the third statement
in the Theorem.

Now we want to describe explicitly the multiplication operation in the
algebra $\scr A$ in terms of the linear tensor product decomposition
(\ref{eq:ATHVH}) of $\scr A$ into the two specified subalgebras of
$\scr A$.  But by the comments before (\ref{eq:CAH}), to prove the
rest all we need to show is that $T(H \cdot V)$ is stable under the
adjoint action of $H$.  But this is clear from the iteration of
(\ref{eq:hcc}) for a product of several elements of $H \cdot V$,
together with (\ref{eq:h1}).  This completes the proof of Theorem 2.1.
$\square$

\begin{remark}

The proof in \cite{BL} proceeds by viewing $T(H \otimes V) \# H$ as a
universal object and by constructing mutually inverse canonical maps
between this structure and $\scr A$.  The argument above essentially
carries this out ``concretely.''

\end{remark}

\end{document}